\documentclass[article,12pt]{amsart}
\usepackage{amsfonts}
\usepackage{amsmath}
\usepackage{graphicx}
\usepackage{color}
\usepackage{epsfig}

\newcommand{\dF}{\mathbb{F}}
\newcommand{\dR}{\mathbb{R}}
\newcommand{\dE}{\mathbb{E}}

\newcommand{\cN}{\mathcal{N}}

\newcommand{\rI}{\mathrm{I}}

\newcommand{\cF}{\mathcal{F}}

\newcommand{\veps}{\varepsilon}

\newcommand{\ind}{\mbox{1}\kern-.25em \mbox{I}}

\def\build#1_#2^#3{\mathrel{\mathop{\kern 0pt#1}\limits_{#2}^{#3}}}

\def\videbox{\mathbin{\vbox{\hrule\hbox{\vrule height1ex \kern.5em
\vrule height1ex}\hrule}}}
\def\demend{\hfill $\videbox$\\}

\numberwithin{equation}{section}
\theoremstyle{plain}

\newtheorem{thm}{Theorem}[section]
\newtheorem{rem}{Remark}[section]
\newtheorem{cor}{Corollary}[section]

\email{bernard.bercu@u-bordeaux.fr}
\email{vvazquez@fcfm.buap.mx}

\begin{document}
\title[Almost sure CLT for ARX processes in adaptive tracking]
{On the almost sure central limit theorem for ARX processes in adaptive tracking}
\author{Bernard Bercu}
	\address{Universit\'e de Bordeaux, Institut de Math\'ematiques, UMR 5251, 351 cours de la lib\'eration, 33405 Talence cedex, France.}

\author{Victor Vazquez}
\address{ Universidad Aut\'onoma de Puebla, Facultad de Ciencias Fisico
Matem\'aticas, Avenida San Claudio y Rio Verde, 72570 Puebla, Mexico.}

\begin{abstract}
The goal of this paper is to highlight the almost sure central limit theorem for martingales to the control community
and to show the usefulness of this result for the system identification of controllable ARX$(p,q)$ process in adaptive tracking. 
We also provide strongly consistent estimators of the even moments of the driven noise of a controllable ARX$(p,q)$ process
as well as quadratic strong laws for the average costs and estimation errors sequences. Our theoretical results are illustrated 
by numerical experiments.
\end{abstract}

\maketitle

\section{INTRODUCTION} \label{intro}

Zadeh \cite{zadeh} introduced the concept of system identification as the determination, on the basis of input and output, 
of a system within a specified class of systems to which the system under test is equivalent. Hence, any theoretical result which leads us to make 
such determination in a more precise way will be a step into the spirit of Zadeh's definition.
In order to track some of the most relevant results in system identification, we may find exhaustive and very useful reviews summarizing the most important contributions in this area
of applied mathematics.
To the best of our knowledge, the most relevant of them are the survey of Astr\"om and Eykhoff \cite{Astrom1971}, the excellent book of Astr\"om and Wittenmark \cite{Astrom},
the work of Ljung \cite{ljung} devoted to adaptive tracking in system identification,
and the beautiful and captivating book-chapter of Gevers \cite{gevers} which deals with the recent developments in identification theory. We also refer the reader to \cite{Hong} for
an overview of basic research on model selection approaches for linear systems and to  \cite{Pillonetto} for Kernel methods in system identification, machine learning and functional estimation.

Since the  pioneer works on system identification and adaptive control, there has been a great deal of activity
from the control community on the theoretical aspects as well as on the practical applications.  Recently, Cho et al. \cite{cho} proposed a 
new parameter estimation method in the framework of composite model reference adaptive control, in order to improve parameter estimation without persistent excitation. 
Heydari \cite{heydari} investigated the stability of adaptive optimal control using value iteration, initiated from a stabilizing control policy. 
Jaramillo et al. \cite{jaramillo} presented an adaptive control framework for compensation of uncertainties and perturbations that satisfy the matching condition on a class of
nonlinear dynamic systems. Moreover, Zhu et al. \cite{zhu} proposed an adaptive model predictive control for unconstrained discrete-time linear systems with parametric uncertainties
We also refer the reader to Gao et al. \cite{gao} who investigated the problem of adaptive tracking control for a class of stochastic uncertain nonlinear systems in the presence of input saturation, 
to Zhao et al. \cite{zhao} who studied the adaptive control for linear systems with set-valued observations in order to track a given periodic target, and to Tao et al. \cite{tao}
for the higher-order tracking properties of model reference adaptive control systems.

Bercu et al.  \cite{BercuVazquez, BV} investigated the asymptotic behavior of the least squares estimator for ARX process in adaptive tracking \cite{Astrom}. 
More precisely, a new notion of strong controllability for multidimensional ARX processes was proposed in \cite{BercuVazquez}. In addition, via a persistently  excited version of the adaptive control, 
it has been shown in \cite{BV} how to avoid this strong controllability condition. In the scalar framework, a serial correlation noise was considered in \cite{BPV1} for ARX processes.
The asymptotic behaviour of the least squares estimator was analyzed together with the almost sure convergence of the Durbin-Watson statistics as well as its asymptotic normality.
It led us to proposed a bilateral statistical test for testing whether or not the serial correlation parameter is equal to some non zero fixed value. Finally, in \cite{BPV2}, the introduction of a persistent excitation in the \"Astrom and Wittenmark adaptive control led us to explore a statistical test for detecting the presence of a serial correlated noise together with the asymptotics of the least squares estimator and of the Durbin-Watson statistics.

The primary goal of this paper is to highlight the almost sure central limit theorem (ASCLT) for martingales to the control community
and to show the usefulness of this result for the system identification of a controllable ARX$(p,q)$ process in adaptive tracking. 
The ASCLT has been widely investigated in stochastic approximation  theory \cite{BCF, Peggy, Pelletier} and in statistics \cite{Bin, Holzmann}.
On the one hand, a large literature is available on the ASCLT for sums of independent random variables 
\cite{Brosamler, LP, Schatte88, Berkes}. On the other hand, it is also possible to find many references on the
ASCLT for martingales \cite{BercuFort, BCF, CM, C, L}. 

Surprisingly, the deep impact of the ASCLT has not deeply reached to the control community. To the best of our knowledge, no reference is available in 
the engineering literature dealing the ASCLT. Hence, the aim of this paper is to show how the 
ASCLT for martingales could provide interesting results for increasing the deepness of the system identification of ARX processes in adaptive tracking. 
The main result in this paper is also related to the estimation error sequence and to the properties of the driven noise, in particular to the almost sure 
estimation of its even moments. Consequently, our result leads us to increase our knowledge on the distribution of the
driven noise of a controllable ARX$(p,q)$ process.

The paper in organized as follows. Section \ref{SCSRX} is devoted to the one-dimensional ARX$(p,q)$ processes in adaptive tracking,
while the ASCLT for the least squares estimator is given in Section \ref{SCMR}. Our theoretical results are illustrated by numerical experiments
in Section \ref{SCNE}. A short conclusion is given in Section \ref{SCC}.
The almost sure central limit theorem for martingales is provided in Appendix A, while all technical proofs are postponed
to Appendix B.

\section{ARX processes}
\label{SCSRX}

In this section, we focus our attention on the one-dimensional ARX$(p,q)$ processes in adaptive tracking, given for all $n\geq 0$ by 
\begin{equation}  
\label{ARX}
A(R)X_{n+1}=B(R)U_{n}+\varepsilon_{n+1}
\end{equation}
where $R$ stands for the shift-back operator, $X_{n}, U_{n}$ and 
$\varepsilon _{n}$ are the system output, input and driven noise,
respectively. The polynomials $A$ and $B$ are given for all $z\in\ \mathbb{C}$
by 
\begin{eqnarray*}
A(z) &=&1-a_{1}z-\cdots -a_{p}z^{p}, \\
B(z) &=&1+b_{1}z+\cdots +b_{q}z^{q},
\end{eqnarray*}
where $a_{i}$ and $b_{j}$ are typically unknown real numbers. 
In all the sequel, we shall make use of the well-known causality assumption on $B$, also known as the minimum phase condition, 
as well as the usual notion of controllability for one-dimensional ARX processes. To be more precise, we assume that the polynomial 
$B(z)$ only has zeros with modulus $>1$ and that polynomials $A(z)-1$ and $B(z)$ are coprime. 
Relation \eqref{ARX} may be rewritten in the compact form
\begin{equation}  
\label{MOD}
X_{n+1}=\theta ^{T}\Phi _{n}+U_{n}+\varepsilon_{n+1}
\end{equation}
where $\theta ^{T}=(a_{1},\ldots ,a_{p},b_{1},\ldots ,b_{q})$ and $\Phi _{n}^T=\left( X_{n}^{p},U_{n-1}^{q}\right)$ 
with $X_{n}^{p}=(X_{n},\ldots ,X_{n-p+1})$ and $U_{n}^{q}=(U_{n},\ldots ,U_{n-q+1})$. We shall assume that the driven noise 
$(\varepsilon _{n})$ is a martingale difference sequence adapted to the filtration $\mathbb{F}=(\mathcal{F}_{n})$
where $\mathcal{F}_{n}$ is for the $\sigma $-algebra of the events occurring up to time $n$, which means that for all $n\geq 0$, 
$\dE[\varepsilon_{n+1} | \cF_n]=0$ a.s. Moreover, we assume that, for all $n\geq 0$, 
$\mathbb{E}[\varepsilon_{n+1}^2|\mathcal{F}_{n}]=\sigma^2 $ a.s. where $\sigma^2>0$. Finally, we assume that $(\varepsilon _{n})$ satisfies, 
for some integer $m\geq 1$ and some real number $a>2m$,
\begin{equation} \label{Hip}
\sup_{n\geq0}\mathbb{E}\bigl[\left|\varepsilon_{n+1}\right|^a|\mathcal{F}_n\bigr]<\infty \hspace{1cm} \text{a.s.}
\end{equation}

The goal of adaptive tracking is to regulate the dynamics of the process $(X_n)$ by forcing the output $X_n$ to track, step by step, 
a predictable reference trajectory $(x_n)$ such that
\begin{equation*}
\sum_{k=1}^n x_k^2=o(n)\hspace{1cm}\text{a.s.}
\end{equation*}
Moreover, at the same time, we shall also estimate the unknown parameter $\theta$. 

First, we focus our attention on the estimation of $\theta$. We shall make use of the least squares estimator which satisfies, for all $n\geq 0$, 
\begin{equation}  
\label{WLS}
\widehat{\theta}_{n+1}=\widehat{\theta}_{n}+S_{n}^{-1}	\Phi_{n}
\Bigl(X_{n+1}-U_{n}-\widehat{\theta}_{n}^{\,T}\Phi_{n}\Bigr),
\end{equation}
$$
S_n= \sum_{k=0}^n \Phi_k \Phi_k^T+I_\delta
$$%
where the initial value $\hat{\theta}_{0}$ may be arbitrarily chosen and $I_\delta$ is the identity matrix of order $\delta=p+q$. 

Next, we are concerned with the choice of the adaptive control sequence $(U_n)$. 
With that aim, we shall make use of the adaptive tracking control proposed by Astr\"om and Wittenmark 
\cite{Astrom} given, for all $n \geq 0$, by 
\begin{equation}  \label{CONTROL1}
U_n = x_{n+1}-\widehat{\theta}_n^{\,T}\,\Phi_n.
\end{equation}
By substituting (\ref{CONTROL1}) into (\ref{MOD}), we obtain the closed-loop
system 
\begin{equation} \label{CLS1}
X_{n+1} - x_{n+1}= \pi_n + \varepsilon_{n+1},
\end{equation}
where prediction error $\pi_n = (\theta - \widehat\theta_n)^{\,T}\Phi_n$. Finally, for each integer $m\geq1$ (this integer will 
be related to the conditional moments of order $2m$ of the driven noise), denote by  $(C_{n}(m))$ and $(G_{n}(m))$  
the sequences of average costs and estimation errors given by
\begin{equation} \label{CNM}
C_{n}(m)=\frac{1}{n}\sum_{k=1}^{n}\bigl(X_{k}-x_{k}\bigr)^{2m}
\end{equation}%
and
\begin{equation} \label{GNM}
G_{n}(m)=\sum_{k=1}^{n}k^{m-1}\bigl\|\widehat\theta_k-\theta\bigr\|^{2m}.
\end{equation}
We assume that the polynomial $B(z)$ only has zeros with modulus $>1$. If $r>1$ is strictly less than the smallest modulus of the
zeros of $B(z)$, then $B(z)$ is invertible in the ball with center
zero and radius $r$ and $B^{-1}(z)$ is a holomorphic function. For all $z\in \mathbb{C}$ such that $|z|\leq r$, denote

\begin{equation}
P(z)=B^{-1}(z)(A(z)-1)=\sum_{k=1}^{\infty }p_{k}z^{k}.  \label{DEFP}
\end{equation}%
All the coefficients $p_k$ may be explicitly calculated as functions of $a_1,\ldots,a_p$ and $b_1,\ldots,b_q$ \cite{BercuVazquez}. 
For any $1\leq i\leq q$, let
\begin{equation}
h_{i}=\sum_{k=i}^\infty p_k p_{k-i+1}.
\end{equation}

In addition, denote by $H$ the square matrix of order $q$,
\begin{equation*}
H=\left( 
\begin{array}{ccccc}
h_{1} & h_{2} & \cdots & h_{q-1} & h_{q} \\ 
h_{2} & h_{1} & h_{2} & \cdots & h_{q-1} \\ 
\cdots & \cdots & \cdots & \cdots & \cdots \\ 
h_{q-1} & \cdots & h_{2} & h_{1} & h_{2} \\ 
h_{q} & h_{q-1} & \cdots & h_{2} & h_{1}%
\end{array}%
\right) .  
\end{equation*}%
Let $K$ be the rectangular matrix of dimension $q\times p$ given, if $p\geq q$, by 

\begin{equation*}  
K=\left( 
\begin{array}{ccccccc}
0 & p_1 & p_2 & \cdots & \cdots & p_{p-2} & p_{p-1} \\ 
0 & 0 & p_1 & \cdots & \cdots & p_{p-3} & p_{p-2} \\ 
\cdots & \cdots & \cdots & \cdots & \cdots & \cdots & \cdots \\ 
0 & \cdots & 0 & p_1 & p_2 & \cdots & p_{p-q+1} \\ 
0 & \cdots & 0 & 0 & p_1 & \cdots & p_{p-q}%
\end{array}
\right)
\end{equation*}
while, if $p\leq q$, by 
\begin{equation*}
K=\left( 
\begin{array}{ccccc}
0 & p_1 & \cdots & p_{p-2} & p_{p-1} \\ 
0 & 0 & p_1 & \cdots & p_{p-2} \\ 
\cdots & \cdots & \cdots & \cdots & \cdots \\ 
0 & \cdots & 0 & 0 &  p_1 \\
\cdots & \cdots & \cdots & \cdots & \cdots \\
0 & 0 & \cdots & 0 & 0 
\end{array}
\right)  .
\end{equation*}

Finally, denote by $L$ the square matrix of order $\delta=p+q$
\begin{equation} 
L=\left( 
\begin{array}{cc} 
I_p & K^{T} \\ 
K & H%
\end{array}%
\right).  \label{DEFl}
\end{equation}

\section{MAIN RESULTS}
\label{SCMR}

Our first result deals with the ASCLT for the least squares estimator $\widehat{\theta}_{n}$ of the unknown parameter $\theta$.

\begin{thm} 
\label{thm1}
Assume that the $\mbox{ARX}(p,q)$ process is causal and controllable. Moreover, assume that for some real number $a>2$,
\begin{equation} \label{Hip3}
\sup_{n\geq0}\mathbb{E}\bigl[\left|\varepsilon_{n+1}\right|^a|\mathcal{F}_n\bigr]<\infty \hspace{1cm} \text{a.s.}
\end{equation}
Then, we have the ASCLT
\begin{equation} \label{ASCLTLS}
\frac{1}{\log n}\sum_{k=1}^n \frac{1}{k}\delta_{\sqrt{k}(\widehat{\theta}_{k}-\theta )} \Rightarrow \cN_\delta(0,L^{-1}) \hspace{1cm} \text{a.s.}
\end{equation}
In other words, for any bounded continuous function $h$,
\begin{equation*}
\vspace{-1ex}
\lim_{n\rightarrow \infty}\frac{1}{\log n}\sum_{k=1}^n \frac{1}{k} h\bigl( \sqrt{k}(\widehat{\theta}_{k}-\theta )\bigr) = \int_{\mathbb{R}^\delta} h(x)dG(x)\hspace{1cm} \text{a.s.}
\end{equation*}
where $G$ stands for the $\cN_\delta(0,L^{-1})$ Gaussian measure. 
\end{thm}

\begin{rem} One can observe that the scalar variance $\sigma^2$ vanishes in the ASCLT.
\end{rem}

\begin{cor} \label{cor1}
Assume that the ARX$(p,q)$ process is causal and controllable. Moreover, assume that $(\varepsilon_n)$ satisfies, for some integer $m\geq1$, condition \eqref{Hip}. 
Then, we have
\begin{equation} 
\label{TH551}
\lim_{n\rightarrow \infty} \frac{1}{\log n}\sum_{k=1}^n k^{m-1} \bigl((\widehat{\theta}_{k}-\theta)^T L (\widehat{\theta}_{k}-\theta)\bigr)^m=\ell(m) \hspace{1cm} \text{a.s.}
\end{equation}
where
\begin{equation}
\label{DEFL2}
\ell(m)= (p+q)\prod_{k=1}^{m-1}(p+q+2k).
\end{equation}
\end{cor}

\begin{rem} \label{rem3}
In the special case $m=1$, we can deduce from \eqref{TH551} the quadratic strong law
\begin{equation} \label{QSL}
\lim_{n\rightarrow \infty} \frac{1}{\log n}\sum_{k=1}^n (\widehat{\theta}_{k}-\theta)^T L (\widehat{\theta}_{k}-\theta)=p+q \hspace{1cm} \text{a.s.}
\end{equation}
\end{rem}

We now focus our attention on the average costs and estimation errors sequences $(C_{n}(m))$ and $(G_{n}(m))$ given by \eqref{CNM} and \eqref{GNM}. 
First of all, it was proven in Lemma $3$ of \cite{BercuVazquez} that the matrix $L$ is positive definite. 
Hence, \eqref{TH551} immediately implies that

\begin{equation*}
G_n(m)=\mathcal{O}(\log n) \hspace{1cm} \text{a.s.}
\end{equation*}
Furthermore, denote 
\begin{equation*}
\Gamma_n(m)=\frac{1}{n}\sum_{k=1}^{n}\varepsilon_{k}^{2m}.
\end{equation*}
The asymptotic behavior of $(C_{n}(m))$ is as follows.

\begin{cor} \label{cor2}
Assume that the ARX$(p,q)$ process is causal and controllable. Moreover, assume that $(\varepsilon_n)$ satisfies, for some integer $m\geq1$, condition \eqref{Hip}. 
Suppose that it exists some integer $1\leq s\leq m$ such that $\mathbb{E}\left[\varepsilon^{2s}_{n+1}|\mathcal{F}_n\right]=\sigma(2s)$ a.s. 
Then, $C_n(s)$ is a strongly consistent estimator of $\sigma(2s)$,

\begin{equation} \label{CONS}
\lim_{n\rightarrow \infty}C_n(s)=\sigma(2s)\hspace{1cm} \text{a.s.}
\end{equation}%
More precisely, for all $0<b<1$ such that $2m<ab$, we have

\begin{equation} \label{CONS1}
\bigl(C_n(s)-\Gamma_n(s)\bigr)^2=o\left(n^{b-1}\right) \hspace{1cm} \text{a.s.}
\end{equation}
\end{cor}

\section{Numerical Experiments}
\label{SCNE}
We provide now some numerical experiments in order to illustrate the most relevant almost sure results of Section \ref{SCMR}. 
More precisely, we shall focus attention on the quadratic  strong law given by \eqref{QSL} as well as on the almost sure convergence of even moments 
given by convergence \eqref{CONS} for different values of $m$. For the sake of simplicity, we assume that the driven noise $(\varepsilon _{n})$ is a sequence 
of independent and identically distributed random variables sharing the same $\mathcal{N}(0,\sigma^2)$ distribution with $\sigma^2=0.8$, 
and the reference trajectory $(x_n)$ is identically zero. Consider the ARX$(2,2)$ process given by \eqref{ARX} where
$$A(z)=1+\frac{6}{5}z-\frac{1}{2}z^2
\hspace{1cm} \text{and}\hspace{1cm}
B(z)=1+\frac{2}{5}z+\frac{1}{4}z^2.$$
One may observe that $B$ is causal since its complex roots have modulus $2$. Moreover, one may easily check that the process is controllable and the matrix 
$L$ is given by

\begin{equation*}
L=
\begin{pmatrix}
1 & 0 & 0 & 0 \\ 
0 & 1 & 6/5 & 0 \\ 
0 & 6/5 & 244/99 & -628/495 \\ 
0 & 0 & -628/495 & 244/99 
\end{pmatrix}.
\end{equation*}
In order to illustrate the quadratic strong law given by \eqref{QSL}, the sample size will increases from $n=100$ to $n=5000$, 
and we shall denote by $\Delta_n$ the average of $N=100$ values of
\begin{equation*}
\frac{1}{\log n}\sum_{k=1}^n (\widehat{\theta}_{k}-\theta)^TL (\widehat{\theta}_{k}-\theta).
\end{equation*}
We may conclude by observing Table $1$ that, even with the slow growth of the logarithmic function, relative errors are small, 
and the quadratic strong law is nicely shown. We recall here that the almost sure limit given by \eqref{QSL} is $p+q=4$.
\medskip

\begin{center}
\begin{tabular}{|c|c|c|}
\hline
 $n$ & $\Delta_n$ & relative error\\
\hline
 $100$ & $4.33$  & $8.14\%$ \\
\hline
 $500$ & $4.14$ & $3.51\%$ \\
\hline
 $1\,000$ & $3.97$ & $0.69\%$ \\
\hline
 $2\,000$ & $4.03$  & $0.78\%$ \\
\hline
$5\,000$ & $3.972$  & $0.70\%$ \\
\hline 
\end{tabular}\\
\vspace{0.5cm}
Table 1. Quadratic strong law.
\end{center}
Let us deal now with the almost sure convergence of even moments given by convergence \eqref{CONS}. For that purpose, we shall consider the average of $N=100$ 
values of sample size $n=10\,000$ of $C_n(m)$. The corresponding results are presented in Table $2$ where the values of $m$ increases from $1$ to $5$.
\medskip

\begin{center}
\begin{tabular}{|c|c|c|c|}
\hline
 $m$ & $C_n(m)$ & $\sigma(2m)$ &relative error\\
\hline
 $1$ & $0.6426$ & $0.64$ & 0.41\%\\
\hline
 $2$ & $1.245$ & $1.229$ & 1.84\%\\
\hline
 $3$ & $4.07$ & $3.932$ & 3.51\%\\
\hline
 $4$ & $18.81$ & $17.62$ & 6.75\% \\
\hline
 $5$ & $110.61$ & $101.47$ & 9.00\% \\
\hline
\end{tabular}\\
\vspace{0.5cm}
Table 2. Convergence of even moments.
\end{center}
We observe that as the value of $m$ increases, the relative error also increases. In other words, it is necessary to take large sample sizes in order 
to estimate large order even moments. For example, choose the value of $m=5$ and consider the large sample sizes 
$n=20\,000$, $30\,000$ and $50\,000$ as indicated in Table $3$. As expected, for large sample sizes values, 
the almost sure convergence of even moments can be improved substantially.
\bigskip

\begin{center}
\begin{tabular}{|c|c|c|c|}
\hline
  n & $C_n(5)$ & relative error\\
\hline
  20\,000 & $108.54$ & 6.90\%\\
\hline
  30\,000 & $106.15$ & 4.61\% \\
\hline
  50\,000 & $105.44$ & 3.91\% \\
\hline
\end{tabular}\\
\vspace{0.5cm}
Table 3. Estimation of $\sigma(10)$.
\end{center}

\section{Conclusion}
\label{SCC}

In this paper, we established the almost sure central limit theorem for the least squares estimator
of the unknown parameter of a controllable ARX$(p,q)$ process in adaptive tracking.  We have also provided a strongly consistent estimator 
for the even moments of the driven noise. Even when most of the engineering methods do not consider high order moments, it is useful to go deeper into 
the knowledge of the driven noise distribution through the estimation of such moments since it gives us a better notion of the underlying uncertainty.

\section*{Appendix A \\ On the almost sure central limit theorem for martingales}
\renewcommand{\thesection}{\Alph{section}}
\renewcommand{\theequation}{\thesection.\arabic{equation}}
\setcounter{section}{1}
\setcounter{equation}{0}

The goal of this Appendix is to highlight the ASCLT for martingales \cite{Bercu2004, BercuFort, BCF, CM, C, L} to the control community.
Let $(\Omega,\cF,\mathbb{P})$ be a probability space endowed with a filtration $\dF= (\cF_{n})$
where $\mathcal{F}_{n}$ is for the $\sigma $-algebra of the events occurring up to time $n$.
Assume that $(M_n)$ be a sequence of integrable random vectors in $\mathbb{R}^d$ such that, for all $n\geq 0$, $M_n$ is $\cF_{n}$-measurable.
We shall say that $(M_n)$ is a martingale with respect to the filtration $\dF$ if for all $n\geq 0$, $\dE[ M_{n+1}|\cF_{n}]=M_n$ almost surely. 
Throughout this Appendix, $(\veps_n)$ is a martingale difference sequence adapted to $\dF$ such that, for all $n \geq 0$,  
$\dE[\veps_{n+1}^{2}|\cF_{n}]=\sigma^2 $ a.s. where $\sigma^2>0$.
Let $(\Phi_{n})$ be a sequence of random vectors of $\dR^d$, adapted
to $\dF$. Denote by $(M_n)$ the locally square integrable martingale
$$
M_{n}=M_0+\sum_{k=1}^{n}\Phi_{k-1}\veps_{k},
$$%
where the initial value $M_0$ can be taken arbitrarily. Its increasing process $\langle M \rangle_n$ is defined, for all $n \geq 1$, by
$$
\langle M \rangle_n = \sum_{k=1}^n \dE[\Delta M_k\Delta M_k^T | \cF_{k-1}]
$$
where $ \Delta M_k=M_k-M_{k-1}$.
We clearly have $\langle M \rangle_n =\sigma^2 S_{n-1}$
where
$$
S_n= \sum_{k=0}^n \Phi_k \Phi_k^T.
$$
A simplified version of the ASCLT for multivariate martingales is as follows \cite{CMT}.
\begin{thm} 
\label{thm2}
Assume that it exists a positive definite symmetric matrix $L$ such that
\begin{equation}
\label{ASCLTC1}
\lim_{n \rightarrow \infty} \frac{1}{n}S_n=L \hspace{1cm} \text{a.s.}
\end{equation}
Moreover, assume that $(M_n)$ satisfies Lindeberg's condition which means that for all
$\veps>0$,
\begin{equation}
\label{Lind}
\lim_{n \rightarrow \infty} \frac{1}{n}
\sum_{k=1}^{n}\dE \bigl[ ||\Delta M_{k}||^{2}\rI_{\left\{||\Delta M_{k}|| \geq \veps \sqrt{n}
\right\}} | \cF_{k-1}\bigr]
=0 \hspace{1cm} \text{a.s.}
\end{equation}
Then, we have the ASCLT
\begin{equation*}
\label{ASCLTM}
\frac{1}{\log n}\sum_{k=1}^n \frac{1}{k}\delta_{M_k/\sqrt{k}} \Rightarrow \cN_d(0,\sigma^2 L) \hspace{1cm} \text{a.s.}
\end{equation*}
In other words, for any bounded continuous function $h$,
\begin{equation*}
\vspace{-1ex}
\lim_{n\rightarrow \infty}\frac{1}{\log n}\sum_{k=1}^n \frac{1}{k} h\left( \frac{M_k}{\sqrt{k}} \right) = \int_{\mathbb{R}^d} h(x)dG(x)\hspace{1cm} \text{a.s.}
\end{equation*}
where $G$ stands for the $\cN_d(0,\sigma^2 L)$ Gaussian measure. 
\end{thm}
The convergence of the even moments in the ASCLT for multivariate martingales was established in \cite{BCF}. 
\begin{thm} 
\label{thm3}
Assume that the almost sure convergence \eqref{ASCLTC1} is satisfied. In addition, 
suppose that for some integer $m\geq 1$ and for some real number $a>2m$,
\begin{equation} 
\label{ASCLTC2}
\sup_{n\geq 0}\dE\bigl[|\veps_{n+1}|^a|\cF_n\bigr]<\infty \hspace{1cm}\text{a.s.}
\end{equation}
Then, we have
\begin{equation}
\label{ASCLTMOM}
\lim_{n\rightarrow \infty}\frac{1}{\log n}\sum_{k=1}^n \frac{1}{k} \bigl(M_k^TS_{k-1}^{-1}M_k\bigr)^m=\ell(m)
\hspace{1cm}\text{a.s.}
\end{equation}
where
\begin{equation}
\label{DEFL}
\ell(m)= d\sigma^{2m}\prod_{k=1}^{m-1}(d+2k).
\end{equation}
\end{thm}

\begin{rem}
The limit $\ell(m)$ corresponds exactly to the mean value of 
$||Z||^{2m}$ where $Z$ has a standard $\cN_d(0,\sigma^2 I_d)$ distribution. 
Consequently, Theorem~\ref{thm3} can be seen as the convergence of moments of order $2m$ 
in the ASCLT for multivariate martingales. 
\end{rem}

\section*{Appendix B \\ Proofs of our main results}
\renewcommand{\thesection}{\Alph{section}}
\renewcommand{\theequation}{\thesection.\arabic{equation}}
\setcounter{section}{2}
\setcounter{equation}{0}

\textit{Proof of Theorem \ref{thm1}.} It follows from \eqref{MOD} and  \eqref{WLS} that for all $n\geq 1$, 
\begin{equation}  
\label{LS22}
\widehat{\theta}_{n}-\theta=S_{n-1}^{-1}M_{n}
\end{equation}
where 
\begin{equation*} 
M_{n}=M_0+\sum_{k=1}^n\Phi_{k-1}\varepsilon_k
\end{equation*}
with $M_0=\widehat{\theta}_0-\theta$. It was proven in Theorem $5$ of \cite{BercuVazquez} that
\begin{equation}  \label{TH11}
\lim_{n\rightarrow \infty} \frac{S_{n}}{n} = \Lambda \hspace{1cm} \text{a.s.}
\end{equation}
where $\Lambda=\sigma^2 L$ and the limiting matrix $L$ is given by \eqref{DEFl}. Moreover, we can deduce from \eqref{Hip3} that Lindeberg's 
condition \eqref{Lind} is satisfied. Consequently, we obtain \eqref{ASCLTLS} from Theorem \ref{thm2} together with \eqref{LS22} and \eqref{TH11}. 
\demend
\newline
\noindent
\textit{Proof of Corollary \ref{cor1}.} The almost sure convergence \eqref{TH551} follows from the conjunction of Theorem \ref{thm1} together with \eqref{LS22} and \eqref{TH11}, using 
the same arguments as in the proof of Corollary $3.3$ in \cite{BCF}. 
\demend
\newline
\noindent
\textit{Proof of Corollary \ref{cor2}.} For any integer $1\leq s \leq m$, we obtain from \eqref{CLS1} that
\begin{eqnarray} \label{eq25}
n\bigl(C_n(s)-\Gamma_n(s)\bigr)&=&\sum_{k=1}^n(X_k-x_k)^{2s}-\sum_{k=1}^n\varepsilon_k^{2s} \nonumber \\
&=&\sum_{k=1}^n(\pi_{k-1}+\varepsilon_k)^{2s}-\sum_{k=1}^n \varepsilon_k^{2s} \nonumber \\
&=&\sum_{k=1}^n \pi_{k-1}^{2s}+R_n(s)
\end{eqnarray}
where
$$
R_n(s)=\sum_{l=1}^{2s-1} \sum_{k=1}^n  \binom{2s}{l} \pi_{k-1}^{2s-l} \varepsilon_k^l.
$$
On the one hand, in the special case $s=1$, we deduce from Theorem \ref{thm3} with $m=1$ that
\begin{equation} \label{eq26}
\lim_{n\rightarrow \infty} \frac{1}{\log n}\sum_{k=1}^n \pi_k^2=\sigma^2(p+q)\hspace{1cm}\text{a.s.}
\end{equation}
In addition, we find from the strong law of large numbers for martingales \cite{Duflo} that
\begin{equation} \label{eq27}
R_n(1)=o(\log n) \hspace{1cm}\text{a.s.}
\end{equation}
Hence, we obtain from \eqref{eq25} together with \eqref{eq26} and \eqref{eq27} that
\begin{equation} \label{eq28}
\lim_{n\rightarrow \infty} \frac{n}{\log n}\bigl(C_n(1)-\Gamma_n(1)\bigr)=\sigma^2(p+q) \hspace{1cm}\text{a.s.}
\end{equation}
Therefore, as
$$\lim_{n\rightarrow \infty} \Gamma_n(1)=\sigma^2 \hspace{1cm}\text{a.s.}$$
convergence \eqref{eq28} clearly leads to \eqref{CONS} since $\sigma(2)=\sigma^2$. On the other hand, for $2\leq s \leq m$, it follows from convergence \eqref{TH11} 
that $\log d_n \sim (p+q)\log n$ where $d_n=\det(S_n)$. Consequently, Corollary $3.1$ of \cite{BCF} lead us to
\begin{equation} \label{eq29}
\lim_{n\rightarrow \infty} \frac{1}{\log n} \sum_{k=1}^n \pi_k^{2s}=0\hspace{1cm}\text{a.s.}
\end{equation}%
Moreover, it is not hard to see that the remainder term $R_n(s)$ plays a negligible role. Hence, as
$$
\lim_{n\rightarrow \infty} \Gamma_n(s)=\sigma(2s) \hspace{1cm}\text{a.s.}
$$
we obtain \eqref{CONS} from \eqref{eq25} and \eqref{eq29}. Finally, we deduce \eqref{CONS1} from Remark $3.1$ 
of \cite{BCF}, which completes the proof of Corollary \ref{cor2}.
\demend

\end{document}